\documentclass[11pt]{article}
\usepackage{amssymb}
\usepackage{amsthm}

\usepackage{amsmath}
\usepackage{amsfonts}
\usepackage{amssymb}

\usepackage{braket}

\newtheorem{Theorem}{Theorem}

\newtheorem{remark}{Remark}

\newcommand{\FF}{\mathbb{F}}
\newcommand{\QED}{\begin{flushright} $\Box$ \end{flushright}}

\begin{document}
\title{Complexity of multivariate polynomial evaluation}
\author{E. Ballico, M. Elia, and M. Sala}
\date{}
\maketitle
    
\begin{abstract}
We describe a method to evaluate multivariate polynomials over a finite field
and discuss its multiplicative complexity.
\end{abstract}



\section{Introduction}
Many applications require the evaluation of multivariate polynomials over finite fields. For instance, the so called affine codes (also called evaluation or functional or algebraic geometry codes) are obtained evaluating a finite-dimensional linear subspace
of $\mathbb {F}_q[x_1,\dots ,x_r]$ at a finite set $S\subseteq \mathbb {F}_q^r$ (\cite{g}, \cite{h}, \cite{l} and several other papers). 
When the degree $n$ of the polynomials is small, and/or the number $r$ of variables is also small the direct computation is
efficient, however as $n$, or $r$, or both become large, evaluation becomes an issue. The case of univariate polynomials was considered by several authors, see e.g \cite{patterson}, \cite{pan} and some recent papers (\cite{ers1}, \cite{ERS}).
In this paper we propose an evaluation method for multivariate polynomials 
which reduces significantly the multiplicative complexity and hence
the computational burden.

\noindent
Set $M_r(n)=\binom{n+r}{r}$, and let $p(x_1, \ldots,x_r)$ be a polynomial of degree $n$ in $r$
 variables with coefficients in a finite field $\mathbb{F}_{p^s}$; the number of monomials occurring in
 $p(x_1, \ldots,x_r)$ is $M_r(n)$. 
We will consider the evaluation of $p(x_1, \ldots,x_r)$ at a point 
 $\mathbf a=(\alpha_1, \ldots, \alpha_r) \in F_{p^m}^{r}$, where $m$ is divisible by $s$.
A direct evaluation of $p(\alpha_1, \ldots, \alpha_r)$ 
 is obtained from the evaluation of the $M_r(n)$ distinct monomials, a task requiring $M_r(n)-r-1$
 multiplications. Therefore we perform $M_r(n)-1$ multiplications, and a total number
 $A_r(n)=M_r(n)-1$ of additions. The total number of required multiplications is
$$  P_r(n)= 2M_r(n)-r-2 = 2\frac{n^r}{r!}+\frac{r+1}{(r-1)!}n^{r-1}+ \cdots -1-r  \,, $$
however different computing strategies can require a significant smaller number of multiplications.
To the aim of developing some of these strategies, the polynomial $p(x_1, \ldots,x_r)$ is written as a sum 
\begin{equation}
   \label{eqint}
   p(x_1, \ldots,x_r) = \sum_{i=0}^{s-1} \beta^i q_i(x_1, \ldots,x_r)
\end{equation}
of $s$ polynomials, where each $q_i(x_1, \ldots,x_r)$ is a polynomial of degree $n$ in $r$ variables
 with coefficients in the prime field $\mathbb {F}_p$, the value $p(\alpha_1, \ldots, \alpha_r)$ can be obtained from the $s$ values
 $q_i(\alpha_1, \ldots, \alpha_r)$, $2s-2$ multiplications, and $s$ additions in $\mathbb {F}_{p^m}$. In these 
 computations $\beta$ and its powers constitute a basis of $\mathbb {F}_{p^m}$.
Therefore, we may restrict our attention to the evaluation at a point $\mathbf a \in \mathbb{F}_{p^m}^{r}$
 of a polynomial $q(x_1, \ldots,x_r)$, in $r$ variables, of degree $n$ over $\mathbb{F}_p$.  

\vspace{3mm}
\noindent
As pointed out in \cite{ERS}, \S 2.1, the prime $2$ is particularly interesting because 
of its occurrence in many practical applications, for example in error correction coding. Furthermore, in
 $\mathbb{F}_2$ multiplications are trivial.  
Therefore, we give first a description of our method in the easiest
case, that is, over $\FF_2$ and with two variables.
Later, we generalize to any setting.

\section{Our computational model}

There are two kinds of multiplications that are involved in our computations:
field multiplications in the coefficient field $\FF_{p^s}$ and in extension
field $\FF_{p^m}$.
We assign cost $1$ to any of these, except for the multiplications
by $0$ or $1$, that cost $0$ in our model.
\begin{remark}
There can be multiplications that cost much less, such as squares
in characteristic $2$, but we still treat them as cost $1$.
\end{remark}
As customary, we assign cost $0$ to any data reading. 

We could count separately field sums, but our aim is to minimize
the number of field multiplications, and so we use as implicit upper
bound for the number of sums the value $2M_r(n)$, that is, twice the number of
all monomials.
We will not discuss of the number of sums any further.

We assume that an ordering on monomials is chosen once and for all,
e.g. the degree lexicographical ordering 
(see \cite{CGC-cd-inbook-D1moratech}), so that
our input data can be modeled as an $\FF_{p^s}$ string, any entry
corresponding to a polynomial coefficient.

\begin{remark} \label{qq}
A well-established method to evaluate all monomials up to degree $n$
at a given point is to start from degree-$1$ monomials and then
iterate from degree-$r$ monomials to degree-$r+1$ monomials,
since the computations of any degree-$r+1$ monomial requires only
one multiplication, once you have in memory all degree-$r$ monomials.
\end{remark}

We remark here that our algorithm accepts as input {\emph any} polynomial
of a given total degree and so our estimates are worst-case complexity,
which translates in considering dense polynomials.
Clearly, other faster methods could be derived for special classes of polynomials, 
such as sparse polynomials or polynomials with a predetermined algebraic structure.

We will not discuss the memory requirement of our methods, but one can
see easily by inspecting the following algorithms that it is negligible
compared to their computational effort.

\section{The case $r=2$, $p=2$}
A polynomial $P(x, y)$ of degree $n$ in $2$ variables over the binary field may be decomposed 
into a sum of $4$ polynomials as
\begin{eqnarray}
 \label{eq0}
P(x, y) &=& P_{0,0} (x^2, y^2)+ xP_{1,0} (x^2, y^2) + yP_{0,1} (x^2, y^2) + x yP_{1,1} (x^2, y^2) \nonumber \\
&=&  P_{0,0} (x, y)^2+ xP_{1,0} (x, y)^2 + yP_{0,1} (x, y)^2 + x yP_{1,1} (x, y)^2.
\end{eqnarray}
where $P_{i, j}(x, y)$ are polynomials of degree $\lfloor \frac{n-i-j}{2}\rfloor$.
Therefore the value of $P(x, y)$ in the point $\mathbf a = (\alpha_1, \alpha_2) \in \mathbb{F}_{2^m}^2$
 can be obtained by computing the $4$ numbers $P_{i, j} (\alpha_1, \alpha_2)$, the monomial 
 $\alpha_1 \alpha_2$, performing $3$ products $\alpha_1 P_{1, 0} (\alpha_1, \alpha _2)$,
 $\alpha_2 P_{0, 1} (\alpha_1, \alpha_2)$, and $\alpha_1 \alpha_2 P_{1, 1} (\alpha_1, \alpha_2)$, 
 and finally performing $3$ additions. 
We observe that all $P_{i,j}$'s have the same possible monomials, i.e. all
monomials of degree up to $\lfloor \frac{n}{2}\rfloor$.
There is no need to store separately $P_{0,0},P_{0,1},P_{1,0},P_{1,1}$,
because the selection of any of these is obtained by a trivial indexing
rule.
 The polynomials $P_{i, j}(\alpha_1, \alpha_2)$ can be evaluated as sums of such monomials,
  which can be evaluated once for all.
Therefore, $P(\alpha_1, \alpha_2)$ is obtained performing (see Remark \ref{qq})
a total number of 
  $$P_r(n)=4+3+ \frac{(\lfloor \frac{n}{2} \rfloor +2)(\lfloor \frac{n}{2} \rfloor +1) }{2} -3 \approx
       \frac{n^2}{8} $$
   multiplications, a figure considerably less than $\frac{n^2}{2}$ as required by the direct computation.
 However, the mechanism can be iterated, and the point is to find the number of steps yielding
  the maximum gain, that is to find the most convenient degree of the polynomials that should be
  directly evaluated. We have the following:
\begin{Theorem}
 \label{theo1}
 Let $P(x,y)$ be a polynomial of degree $n$ over $\mathbb{F}_2$, its evaluation at a point 
  $(\alpha_1, \alpha_2) \in \mathbb{F}_{2^m}^2$ performed by applying repeatedly the
  decomposition (\ref{eq0}),  requires a number $G_2(n,2,L_{opt})$ of products which asymptotically is
$$
G_2(n,2,L_{opt}) \approx c \sqrt{\frac{7}{6} } n~~~~c < 5 ~~~~~~~~~~.
$$
where $L_{opt}$, the number of iterations yielding the minimum of $G_2(n,2,L)$,
 is an integer included into the interval
$$     -  \frac{1}{2} + \log_4 (\sqrt{\frac{6}{7}} n) +\epsilon < {L_{op}} < 
          \log_4 ( \sqrt{\frac{6}{7}}) n +\epsilon' ~~, $$  
where $\epsilon$ and $\epsilon'$ are less than $1$ and $O(\frac{1}{\sqrt n})$.
\end{Theorem}

\noindent
{\sc Proof}. The polynomial $P(x,y)$ is decomposed into the sum of $4$ polynomials that are
perfect squares over $\mathbb{F}_2$, each of which is the similarly decomposed.
 Let $P_{i, j}^{(L,h)} (x, y)$ denote the polynomials at the $L$-step of this iterative process,
 with $h$ varying from $1$ to $4^{L-1}$. The number of polynomials after $L$ steps is $4^L$,
 while their degrees are not greater than $\lfloor \frac{n}{2^L}\rfloor$.
The value $P(\alpha_1, \alpha_2)$ is obtained performing backward the reconstruction process
 obtaining at each step the values $P_{i,j}^{(\ell-1,h)} (\alpha_1, \alpha_2)$ from the values
 $P_{i,j}^{(\ell,h)} (\alpha_1, \alpha_2)$, whereas the $4^L$ numbers $P_{i,j}^{(L,h)} (\alpha_1, \alpha_2)$, $i,j \in \{0,1 \}$ and $h=0, \ldots , 4^{L-1}$,  are computed from the direct
 evaluation of $M_2(\lfloor \frac{n}{2^L}\rfloor )$ monomials using $M_2(\lfloor \frac{n}{2^L}\rfloor )-3$ multiplications.

 Therefore the total number of multiplications necessary to obtain $P(\alpha_1, \alpha_2)$ is a sum of
 $M_2(\lfloor \frac{n}{2^L}\rfloor)-3$ with
\begin{itemize}
  \item[-] the number of squares
$$  \frac{4}{3}(4^{L}-1) = [4^{L}+ 4^{L-1}+\cdots +4^{L-L+1} ]    $$
  \item[-] the number of multiplications of kind $x^iy^j P_{i,j}(\alpha_1, \alpha_2)$
$$ 4^{L}-1= 3 [4^{L-1}+ 4^{L-2}+\cdots +4^{L-L} ]   $$
\end{itemize}
 that is the total number is:
$$ G_2(n,2,L) = \frac{7}{3}(4^{L}-1)+\frac{(\lfloor \frac{n}{2^L} \rfloor +1)(\lfloor \frac{n}{2^L} \rfloor +2)}{2}-3  $$
%
The number of products required to evaluate $P(\alpha_1, \alpha_2)$ in this way is a function
 of $L$, and the values of $L$ that correspond to local minima are specified by the conditions
$$ G_2(n,2,L) \leq G_2(n,2, L-1) ~~\mbox{and}~~ G_2(n,2,L) \leq G_2(n,2,L + 1)   ~~, $$
from which, it is straightforward to obtain the conditions
$$ \begin{array}{l}
  4^L -  \frac{3}{56} \frac{n^2}{4^L} > \frac{n}{2^L14}( \frac{3}{2}-\{\frac{n}{2^L}\})+ (\{\frac{n}{2^L}\}- \{\frac{n}{2^L 2}\}) \left(\frac{-3}{14}+ \frac{1}{14} (\{\frac{n}{2^L}\}+\{\frac{n}{2^L 2}\})+ \frac{1}{2} \right)       \\
            \\
     4^L -  \frac{6}{7} \frac{n^2}{4^L} <
     \frac{4n}{2^L7} (\frac{3}{2}+\{\frac{n}{2^L}\}-2 \{\frac{2n}{2^L}\})- \frac{2}{7} (\{\frac{n}{2^L}\}^2-
      \{\frac{2n}{2^L}\}^2) +\frac{6}{7} (\{\frac{n}{2^L}\}- \{\frac{2n}{2^L}\}) \\
   \end{array} $$ 
where $\{x\}$ denotes the fractional part of $x$. 
These inequalities show that there is only one minimum that corresponds to a value of $L$ such that
$$     -  \frac{1}{2} + \log_4 (\sqrt{\frac{6}{7}} n) +\epsilon < {L_{op}} < 
          \log_4 ( \sqrt{\frac{6}{7}}) n +\epsilon' ~~, $$  
where $\epsilon$ and $\epsilon'$ are $O(\frac{1}{\sqrt n})$.
Therefore, the minimum value of $G_2(n,2,L)$ is asymptotically 
$$  G_2(n,2,L_{op}) \approx c \sqrt{\frac{7}{6}} n  $$
where $c$ is a constant less than $5$. \QED

\begin{remark}\label{b1}
In the computations of our bounds we essentially compute separately each monomial. Hence our approach seems to be very efficient for the computation of several polynomials at the same point. This fact is exploited
 in the computation of the required number of multiplications when the polynomial coefficients are in $\mathbb{F}_{2^s}$. An application of equation (\ref{eqint}) and Theorem \ref{theo1} would give
 the asymptotic estimate
$$  G_{2^s}(n,2,L_{op}) \approx c \sqrt{\frac{7}{6}} n  s~~. $$
since the evaluation of any $q_i$ would cost $c \sqrt{\frac{7}{6}}n$.
 However, the polynomials $q_i(x,y)$ can be evaluated contemporarily. Therefore, computing the power necessary
  to evaluate the polynomial at step $L$ only once, this lead to a total number or required multiplications
$$ G_{2^s}(n,2,L) = s\frac{7}{3}(4^{L}-1)+\frac{(\lfloor \frac{n}{2^L} \rfloor+1)(\lfloor \frac{n}{2^L} \rfloor+2)}{2}-3  $$
because only the reconstruction operations need to be repeated $s$ times. By repeating the argument outlined in the proof of Theorem \ref{theo1}, the conclusion is that the optimal value of $L$ depends also on $s$ and  asymptotically the required value of multiplications is
$$  G_{2^s}(n,2,L_{op}) \approx c' \sqrt{\frac{7}{6}} n \sqrt s ~~. $$
\end{remark}

 
\section{The case $r>2$, $p=2$}

The evaluation of a polynomial $P(x_1, \ldots, x_r)$ in $r$ variables can be done writing $P$, similarly to equation (\ref{eq0}), 
this polynomial as a sum of $2^r$  polynomials 
\begin{eqnarray}
  \label{eq1}
P(x_1, \ldots, x_r) &=& \sum_{i_1, \ldots, i_r \in \{0,1 \}}
x_1^{i_1} \ldots x_r^{i_r} P_{i_1, \ldots, i_r} (x_1^{2}, \ldots, x_r^{2})  \nonumber \\
&=&  \sum_{i_1, \ldots, i_r  \in  \{0,1\}} x_1^{i_1} \ldots x_r^{i_r} 
\left( P_{i_1, \ldots, i_r} (x_1, \ldots, x_r) \right)^2\,,
\end{eqnarray}
where $P_{i_1, \ldots, i_r} (x_1, \ldots, x_r)$ is a polynomial of degree 
$\frac{n-\sum i_j}{2}$.
The argument of Theorem \ref{theo1} still applies, and the minimum number of steps  
 is obtained in the following theorem.
\begin{Theorem}
 \label{theo2}
 Let $L_{opt}$ be the number of steps of this method yielding the minimum number of products, $G_2(n,r,L_{op})$, required to evaluate a polynomial of degree $n$ in $r$ variables,
  with coefficients in $\mathbb{F}_2$. Then $L_{opt}$ is an integer that asymptotically is
  included into the interval
$$ -\frac{1}{2} + \frac{1}{2r} \log_2 \frac{2^r-1}{r!(2^{r+1}-1)} +\frac{\log_2 n}{2}  \leq  L_{op} \leq 
   \frac{1}{2} + \frac{1}{2r} \log_2 \frac{2^r-1}{r!(2^{r+1}-1)} +\frac{\log_2 n}{2}   $$
that is $L_{op}$ is the integer closest to $\frac{1}{2r} \log_2 \frac{2^r-1}{r!(2^{r+1}-1)} +\frac{\log_2 n}{2}$. 
Asymptotically the minimum $G_2(n,r,L_{op})$ is included into the interval:
$$ \frac{1}{\sqrt{2^r}} \sqrt{4\frac{2^{r+1}-1}{2^r-1} \frac{1}{r!}}  n^{r/2} <  G_2(n,r,L_{op}) <
            \sqrt{2^r}  \sqrt{4\frac{2^{r+1}-1}{2^r-1} \frac{1}{r!}}  n^{r/2}~~. $$
\end{Theorem}
\noindent
{\sc Proof}. Using equation (\ref{eq1}) the polynomial $P(x_1, \ldots, x_r)$ evaluated at the point
 $\mathbf a = (\alpha_1, \ldots, \alpha_r) \in \mathbb{F}_{2^m}^r$
 can be obtained from the evaluation of all $P_{i_1, \ldots, i_r} (x_1, \ldots, x_r)$ at $\mathbf a$,
 by  evaluating  $2^r$ monomials  $\alpha_1^{i_1} \ldots \alpha_r^{i_r}$ (which
 require $2^r-r-1$ multiplications),  performing $2^r$ squaring, 
 combining these factors with $2^r-1$ multiplications, and finally
adding the results. \\
We can iterate this procedure: at each step the number of polynomials 
$P_{i,j}$'s is multiplied by $2^r$ and their
degrees are divided at least by $2$. Therefore, after $L$ steps the number of polynomials is $2^{Lr}$ and their degrees are not greater than $\lfloor \frac{n}{2^L} \rfloor$. Once the $2^r$ numbers $P_{i_1, \ldots, i_r} (\alpha_1, \ldots, \alpha_r)$ are known, the total number of squaring is
$$   2^{rL}+ 2^{r(L-1)}+\cdots +2^{r(L-L+1))} = \frac{2^r}{2^r-1} (2^{rL}-1) $$
and  the number of products necessary to obtain $P(\alpha)$ is 
$$  (2^r-1) [2^{r(L-1)}+ 2^{r(L-2)}+\cdots +2^{r(L-L))} ]=2^{rL}-1 \,, $$
hence the total number of required multiplications is
$$ \frac{2^{r+1}-1}{2^r-1} (2^{rL}-1)\,. $$ 
Since the total number of monomials in $r$ variables in a generic polynomial of degree
 $\lfloor \frac{n}{2^L} \rfloor$ is $M_r(\lfloor \frac{n}{2^L} \rfloor)$,
 then $M_r(\lfloor \frac{n}{2^L} \rfloor)-r-1$ is the number of products necessary to evaluate
 all independent monomials.
Therefore, the total number of multiplications for evaluating $P(\mathbf a)$ is
$$  G_2(n,r,L) = \frac{2^{r+1}-1}{2^r-1} (2^{rL}-1)+M_r(n)-r-1 ~~.  $$
We look for the optimal value $L_{op}$ giving the minimum $ G_2(n,r,L_{op})$.
Since
$$  M_r(\lfloor \frac{n}{2^L} \rfloor)=\frac{1}{r!} \left( \frac{n}{2^L}-\{ \frac{n}{2^L} \} \right)^r \prod_{j=1}^{r} \left( 1+ \frac{j}{\frac{n}{2^L}-\{ \frac{n}{2^L} \} } \right)  ~~,  $$
then $M_r(\lfloor \frac{n}{2^L} \rfloor)$ is an expression that is $\frac{1}{r!}(\frac{n}{2^L})^r+ O(\frac{2^L}{n})$ asymptotically in $n$. \\
The local optima are given by the values of $L$ such  the 
$$ G_2(n,r,L) \leq G_2(n,r, L-1) ~~\mbox{and}~~ G_2(n,r,L) \leq G_2(n,r,L + 1)   \ . $$
Then, considering the asymptotic expression  
$$ G_2(n,r,L) = \frac{2^{r+1}-1}{2^r-1} 2^{rL} + \frac{1}{r!} (\frac{n}{2^L})^r ~~, $$ 
it is immediate 
 to obtain the conditions
$$ \begin{array}{l}
 \displaystyle 2^{2rL} > \frac{1}{r!}  \frac{n^r}{2^r} \frac{2^r-1}{2^{r+1}-1}   \\
            \\
 \displaystyle    2^{2rL} < 2^r n^r \frac{1}{r!} \frac{2^r-1}{2^{r+1}-1} ~~,  \\
   \end{array} $$ 
showing that asymptotically $L_{op}$ must satisfy the inequalities
$$ -\frac{1}{2} + \frac{1}{2r} \log_2 \frac{2^r-1}{r!(2^{r+1}-1)} +\frac{\log_2 n}{2}  < L_{op} < 
   \frac{1}{2} + \frac{1}{2r} \log_2 \frac{2^r-1}{r!(2^{r+1}-1)} +\frac{\log_2 n}{2} ~~.  $$
Therefore $L_{op}$ is the closest integer to 
  $\frac{1}{2r} \log_2 \frac{2^r-1}{r!(2^{r+1}-1)} +\frac{1}{2} \log_2 n  $, 
and the total number of products asymptotically is included into the interval:
$$ \frac{1}{\sqrt{2^r}} \sqrt{4\frac{2^{r+1}-1}{2^r-1} \frac{1}{r!}}  n^{r/2} <  G_2(n,r,L_{op}) <  \sqrt{2^r}  \sqrt{4\frac{2^{r+1}-1}{2^r-1} \frac{1}{r!}}  n^{r/2}~~. $$
\QED


\section{The case $r\geq 2$, $p>2$}

A polynomial $P(x_1, \ldots, x_r)$ of degree $n$, in $r$ variables over the field $\mathbb{F}_p$,
 is simply decomposed into a sum of $p^r$ polynomials as
\begin{eqnarray}
  \label{eq2}
P(x_1, \ldots, x_r) &=& \sum_{i_1, \ldots, i_r \in \{0,1, \ldots, p-1 \}}
x_1^{i_1} \ldots x_r^{i_r} P_{i_1, \ldots, i_r} (x_1^{p}, \ldots, x_r^{p})  \nonumber \\
&=&  \sum_{i_1, \ldots, i_r  \in  \{0,1, \ldots, p-1\}} x_1^{i_1} \ldots x_r^{i_r} 
\left( P_{i_1, \ldots, i_r} (x_1, \ldots, x_r) \right)^p
\end{eqnarray}
where $P_{i_1, \ldots, i_r} (x_1, \ldots, x_r)$ is a polynomial of degree 
 $\lfloor \frac{n-\sum i_j}{p}\rfloor$.
Therefore the polynomial $P(x_1, \ldots, x_r)$ evaluated at the point
 $\mathbf a = (\alpha_1, \ldots, \alpha_r) \in \mathbb{F}_{p^m}^r$
can be obtained from the evaluation of all polynomials $P_{i_1, \ldots, i_r} (x_1, \ldots, x_r)$
 at $\mathbf a$, by  evaluating the $p^r$ monomials  $\alpha_1^{i_1} \ldots \alpha_r^{i_r}$ (which require $p^r-r-1$ multiplications),  performing $p^r$ computations of $p$-powers, 
 combining these factors with $p^r$ multiplications, and finally adding all results.\\
The argument of Theorem \ref{theo1} and \ref{theo2} still applies, and the minimum number
 of steps is obtained in the following theorem.
\begin{Theorem}
 \label{theo3}
 Let $L_{opt}$ be the number of steps of this method yielding the minimum number of products, $G_p(n,r,L_{op})$, required to evaluate a polynomial of degree $n$ in $r$ variables,
  with coefficients in $\mathbb{F}_p$. Then $L_{opt}$ is an integer that asymptotically is
  included into the interval
$$ -\frac{1}{2} +B +\frac{\log_p n}{2} \leq L_{op}
 \leq \frac{1}{2} + B +\frac{\log_p n}{2}   $$
where $B=\frac{1}{2r} \log_p \frac{(p-1)(p^r-1)}{r!(2~ p^{r}-1)}$,
that is, $L_{op}$ is the integer closest to $B +\frac{\log_p n}{2}$. 
Asymptotically the minimum $G_p(n,r,L_{op})$ is included into the interval:
$$ \frac{2}{\sqrt{p^r}} \sqrt{(p-1)\frac{2p^{r}-1}{p^r-1} \frac{1}{r!}}  n^{r/2} <  G_p(n,r,L_{op}) <
            2\sqrt{p^r}  \sqrt{(p-1)\frac{2p^{r}-1}{p^r-1} \frac{1}{r!}}  n^{r/2}~~. $$
\end{Theorem}
\noindent
{\sc Proof}. Using equation (\ref{eq2}) the polynomial $P(x_1, \ldots, x_r)$ evaluated at the point
 $\mathbf a = (\alpha_1, \ldots, \alpha_r) \in \mathbb{F}_{p^m}^r$
 can be obtained from the evaluation of all $P_{i_1, \ldots, i_r} (x_1, \ldots, x_r)$ at $\mathbf a$,
 by  evaluating  $p^r$ monomials  $\alpha_1^{i_1} \ldots \alpha_r^{i_r}$ (which
 require $p^r-r-1$ multiplications), computing $p^r$ $p$-powers,
 combining these factors with $p^r-1$ multiplications, and finally performing
the required additions. \\
We can iterate this procedure: at each step the number of polynomials is multiplied by $p^r$ and their
degrees are at least divided by $p$. Therefore, after $L$ steps the number of polynomials is $p^{rL}$ and their degrees are not greater than $\lfloor \frac{n}{p^L} \rfloor$. Once the $p^r$ numbers $P_{i_1, \ldots, i_r} (\alpha_1, \ldots, \alpha_r)$ are known, the total number of $p$-powers is
$$ p^{rL}+ p^{r(L-1)}+\cdots +p^{r(L-L+1))} =\frac{p^r}{p^r-1} (p^{rL}-1)  $$
and  the number of products necessary to obtain $P(\alpha)$ is 
$$ (p^r-1) [p^{r(L-1)}+ p^{r(L-2)}+\cdots +p^{r(L-L))} ]=p^{rL}-1 ~~, $$
hence the total number of required multiplications is
$$ \frac{2p^{r}-1}{p^r-1} (p^{rL}-1)~~. $$ 
The total number of multiplications for computing all the monomials of all the polynomials arising at step $L$ is $M_r(\lfloor \frac{n}{2^L} \rfloor)-r-1$, and further $(p-2) (M_r(\lfloor \frac{n}{2^L} \rfloor)-r-1) $ products are necessary to provide every possible term
 occurring in the polynomials at step $L$. As a consequence the total number of multiplications necessary to evaluate
 $P(\mathbf a)$ is   
$$   \frac{2p^{r}-1}{p^r-1} (p^{rL}-1) + (p-1) (M_r(\lfloor \frac{n}{2^L} \rfloor)-r-1) \,.$$
The same passages used in Theorem \ref{theo2} allow us to conclude that $L_{op}$ is asymptotically
 identified by the chain of inequalities
$$ \sqrt{\frac{1}{p^r}} \sqrt{\frac{(p-1)(p^r-1)n^r}{r!(2p^r-1)}} \leq p^{rL_{op}}  \leq \sqrt{\frac{(p-1)(p^r-1)n^r}{r!(2p^r-1)}} \sqrt{p^r} $$ 
which written in the form  
$$ -\frac{1}{2} +B +\frac{\log_p n}{2} \leq L_{op}
 \leq \frac{1}{2} + B +\frac{\log_p n}{2}   $$
shows that the unique optimal value is the integer closest  to $ B +\frac{\log_p n}{2} $, where $B=\frac{1}{2r} \log_p \frac{(p-1)(p^r-1)}{r!(2~ p^{r}-1)}$.
The minimum number of multiplications is asymptotically included into the interval
$$ \frac{2}{\sqrt{p^r}} \sqrt{(p-1)\frac{2p^{r}-1}{p^r-1} \frac{1}{r!}}  n^{r/2} <  G_p(n,r,L_{op}) <
            2\sqrt{p^r}  \sqrt{(p-1)\frac{2p^{r}-1}{p^r-1} \frac{1}{r!}}  n^{r/2}~~. $$

\begin{remark}\label{b2}
Our proofs start with the evaluations of certain monomials. Hence they may be extended verbatim to other finite-dimensional linear subspaces
of $\mathbb {F}_{p}[x_1,\dots ,x_r]$, just taking their dimension $\alpha$ as vector spaces instead of the integer $\binom{n+r}{r}$. For a suitable
linear space $V$ in Theorem 3 we could get a bound of order
$c_3\sqrt{\alpha}$ with $c_3\sim 2\sqrt{2p^{r+1}}$. For instance, call
$V(r,n)$ the linear subspace of $\mathbb {F}_{p}[x_1,\dots ,x_r]$ formed by all polynomials whose degree in each variable is at most $n$.
We have $\dim (V(n,r)) = (n+1)^r$. In this case iterating this procedure we arrive at each step at a vector
space $V(\lceil n/p^L\rceil ,r)$. Taking $L$ such that $2p^{rL} \sim p\dim (V(\lceil n/p^L\rceil,r))$, i.e. taking $L\sim \log _p n/2 + B$
with $B\sim \frac{1}{2r}\log _p(p)/2 \sim \frac{1}{2r} \sim - \frac{1}{2r}\log _p 2$, we get an upper bound
of order $2\sqrt{2p^{r+1}}n^{r/2}$.
\end{remark}

\section{Further remarks}

The complexity of polynomial evaluation is crucial in the determination
of the complexity of several computational algebra methods, such as
the Buchberger-Moeller algorithm 
(\cite{CGC-alg-art-buchmoeller,CGC-cd-inbook-D1morafglm}), 
other commutative algebra methods (\cite{CGC-cd-inbook-D1moratech}),
the Berlekamp-Massey-Sakata algorithm (\cite{CGC-cd-art-sakata88,CGC-cd-inbook-D1sakata1}).

In turn, these algorithms are the main tools used in algebraic
coding theory (and in cryptography).
This justifies our special interest in the finite field case.
For example, the previous algorithms can be adapted naturally
to achieve  iterative decoding of algebraic codes and algebraic-geometry
codes, see e.g. \cite{CGC-cd-inbook-D1sakata2,CGC-cd-inbook-D1eleanna}.
Other versions can decode and construct more general geometric
codes, see e.g. \cite{g}.

\section{Acknowledgments}
The authors would like to thank C. Fontanari for valuable discussions.

The first author would like to thank MIUR and GNSAGA of INDAM.
The third author would like to thank MIUR for the programme
``Rientro dei cervelli''.

This work has been partially done while the second author was Visiting Professor
with the University of Trento, funded by CIRM.

\vspace{0.5cm}

\noindent
Edoardo Ballico, Massimiliano Sala \newline
Dipartimento di Matematica, Universit\`a degli Studi di Trento \newline 
Via Sommarive 14, 38123 Trento, Italy. \newline
E-mail addresses: ballico@science.unitn.it, \newline
maxsalacodes@gmail.com 

\vspace{0.2cm}

\noindent
Michele Elia \\
Dipartimento di Elettronica, Politecnico di Torino \newline
Corso Duca degli Abruzzi 24, 10129 Torino, Italy \newline
E-mail address: elia@polito.it

\end{document}